\documentclass[12pt]{article}

\usepackage{amssymb}

\newtheorem {theorem} {Theorem}
\newtheorem {definition} [theorem]{Definition}

\newcommand{\bbox}{\ \hfill\rule[-1mm]{2mm}{3.2mm}}

\def\div{{\rm div}}

\title{Linearizable ordinary differential equations. \thanks{The authors are partially
supported by a MCYT/FEDER grant number MTM2005-06098-C02-02. The
second author is also partially supported by a CIRIT grant number
2005SGR 00550, and by DIUE of Government of Catalonia
``Distinci\'o de la Generalitat de Catalunya per a la promoci\'o
de la recerca universit\`aria".}}

\author{{\sc H\'ector Giacomini$^{\ (1)}$, Jaume Gin\'e$^{\ (2)}$ and Maite Grau$^{\ (2)}$}}
\date{}

\begin{document}
\maketitle

\begin{abstract}
Our purpose in this paper is to study when a planar differential
system polynomial in one variable linearizes in the sense that it
has an inverse integrating factor which can be constructed by
means of the solutions of linear differential equations. We give
several families of differential systems which illustrate how the
integrability of the system passes through the solutions of a
linear differential equation. At the end of the work, we describe
some families of differential systems which are Darboux integrable
and whose inverse integrating factor is constructed using the
solutions of a second--order linear differential equation defining
a family of orthogonal polynomials.
\end{abstract}

{\small{\noindent 2000 {\it AMS Subject Classification:} 14H05, 34A05, 34A34.  }}

\section{Introduction \label{sect1}}

In this work we consider planar polynomial differential systems
as:
\begin{equation}
\dot{x}= P(x,y), \quad  \dot{y}=Q(x,y), \label{eq1}
\end{equation}
where $P(x,y)$ and $Q(x,y)$ belong to the ring of real polynomials
in one variable and are analytic in the other variable, that is,
they belong to the ring $\mathbb{R}(x)[y]$ if we choose $y$ as the
variable in which they are polynomial. We will always assume that
$P(x,y)$ and $Q(x,y)$ are coprime polynomials with respect to $y$.
We denote by ${\rm d}$ the maximum of the degrees of $P$ and $Q$
as polynomials in $y$.
\par
We define the {\em orbital equation} associated to system
(\ref{eq1}) as the ordinary differential equation which is
satisfied by the orbits of the system, that is, the orbital
equation associated to system (\ref{eq1}) is either \[
\frac{dy}{dx} \, = \, \frac{Q(x,y)}{P(x,y)}\, , \qquad \mbox{or}
\qquad \frac{dx}{dy} \, = \, \frac{P(x,y)}{Q(x,y)}\, . \] \par The
aim of this work is to study when a system (\ref{eq1}) linearizes.
\begin{definition}
We say that system {\rm (\ref{eq1})} {\em linearizes}, or that it
is linearizable, if it has an inverse integrating factor which can
be constructed by means of the solutions of linear differential
equations. \label{def1}
\end{definition}

We recall that the classical definition that a system (\ref{eq1})
is linearizable is that there exists a change of variables which
transforms the orbital equation associated to system (\ref{eq1})
into a linear differential equation. The techniques used to find
such a change usually come from the Lie group theory, see
\cite{BlumanKumei, Olver} and the references therein. We do not
treat this problem in this paper but the examples that we study
show that there is a connection between both definitions of
linearizability. \newline

This paper is related to the integrability problem which is
defined as the problem of finding a first integral for a planar
differential system and determining the functional class it must
belong to. We recall that a first integral $H(x,y)$ of system
(\ref{eq1}) is a function of class $\mathcal{C}^1$ in some open
set $\mathcal{U}$ of $\mathbb{R}^2$, non locally constant and
which satisfies the following partial differential equation:
\[ P(x,y) \frac{\partial H}{\partial x}(x,y) \, + \, Q(x,y)
\frac{\partial H}{\partial y}(x,y) \, \equiv \, 0. \] An inverse
integrating factor of system (\ref{eq1}) is a function $V(x,y)$ of
class $\mathcal{C}^1$ in some open set $\mathcal{U}$ of
$\mathbb{R}^2$, non locally null and which satisfies the following
partial differential equation:
\[ P(x,y) \frac{\partial V}{\partial x}(x,y) + Q(x,y) \frac{\partial
V}{\partial y}(x,y) = \left( \frac{\partial P}{\partial x}(x,y) +
\frac{\partial Q}{\partial y}(x,y) \right) \,V(x,y) . \] The
function $(\partial P/\partial x) \, + \, (\partial Q/\partial y)$
is called the {\em divergence} of system (\ref{eq1}) and it is
denoted by $\div$ throughout the rest of the paper. We note that
the function $1/V(x,y)$ is an integrating factor for system
(\ref{eq1}) in $\mathcal{U}$, and that given an inverse
integrating factor defined in $\mathcal{U}$, a first integral in
$\mathcal{U}-\{V=0\}$ can be constructed by means of the following
line integral: \[ H(x,y) \, = \, \int_{(x_0,y_0)}^{(x,y)}
\frac{Q(x,y) dx - P(x,y) dy}{V(x,y)}, \] where $(x_0,y_0)$ is any
chosen base point in $\mathcal{U}$ with $V(x_0,y_0) \neq 0$. We
note that this function $H(x,y)$ is only well--defined, in
general, in a simply--connected subset of $\mathcal{U}-\{V=0\}$.
In nonsimply--connected subsets of $\mathcal{U}-\{V=0\}$, $H(x,y)$
can be a multivalued function but it continues to exhibit the
dynamic behavior of the orbits in the set. \par The integrability
of system (\ref{eq1}) is given, in many occasions, by the
existence of invariant curves. We say that a $\mathcal{C}^1$
function $f: \mathcal{U} \subset \mathbb{R}^2 \to \mathbb{R}$ is
an {\em invariant curve} for a system (\ref{eq1}) if it is not
locally constant and satisfies \[ P(x,y) \frac{\partial
f}{\partial x}(x,y) + Q(x,y) \frac{\partial f}{\partial y}(x,y) =
k_f(x,y) \, f(x,y),
\] with $k_f(x,y)$ a polynomial in $y$ of degree lower or equal
than ${\rm d}-1$, where {\rm d} is the degree of the system in
$y$, and it is of class $\mathcal{C}^1$ in the other variable.
This function $k_f(x,y)$ is called the {\em cofactor} of $f(x,y)$.
In case $f(x,y)=0$ defines a curve in the real plane, this
definition implies that the vector field associated to system
(\ref{eq1}) is tangent to the curve $f=0$. In case $f(x,y)$ is a
polynomial we say that $f(x,y)=0$ is an {\em invariant algebraic
curve} for system (\ref{eq1}). The construction of inverse
integrating factors or analytic first integrals inside certain
functional classes (polynomial, rational, elementary or
Liouvillian) is strongly related with the existence of invariant
algebraic curves, see for instance the work \cite{GGG} and
specially the references therein, and it belongs to the context of
the Darboux theory of integrability. \par When considering the
integrability problem we are also addressed to study how the
existence of a first integral in a certain functional class
implies the existence of an inverse integrating factor inside a
certain given class of functions. In the particular case that
system (\ref{eq1}) is polynomial, we have that the existence of an
elementary first integral implies the existence of an inverse
integrating factor which is a rational function up to a rational
power. Moreover, when system (\ref{eq1}) is a polynomial, we have
that the existence of a Liouvillian first integral implies the
existence of an inverse integrating factor of Darboux type, see
\cite{ChGGLl, GGG} and the references therein for the proof of
these results. A function of the form $\exp{(h/g)} \,
f_1^{\lambda_1}f_2^{\lambda_2} \ldots f_s^{\lambda_s}, $ where
$h(x,y)$ and $g(x,y)$ are polynomials, each $f_i(x,y)=0$ is an
invariant algebraic curve for system (\ref{eq1}), $\lambda_i \in
\mathbb{C}$ not all of them null, for $i=1,2,\ldots,s$ and $s \in
\mathbb{N}$, is called a {\em Darboux} function. These results
suggest that the functional class of an inverse integrating factor
is usually easier than the functional class of a first integral.
This is the reason why we look for an inverse integrating factor
to study the integrability of system (\ref{eq1}). Moreover, the
inverse integrating factor is shown to be defined in phase
portraits in which the dynamics avoid the existence of a first
integral. However, there are also systems whose dynamics avoid the
existence of an inverse integrating factor, see \cite{flows}.
\newline

In the work \cite{GGG}, systems of the form (\ref{eq1}) whose
integrability is given by the solutions of linear differential
equations are described. We obtained a result which allows to find
an explicit expression for a first integral of a certain type. By
means of a rational change of variable, we made correspond the
homogenous second order linear differential equation:  $A_2 (x)
w''(x) + A_1(x)  w'(x) + A_0(x) w(x) = 0$, whose coefficients are
polynomials, to a planar polynomial differential system. We prove
that this system has an invariant curve for each arbitrary non
null solution $w(x)$ of the second-order ordinary differential
equation, which, in case $w(x)$ is a polynomial, gives rise to an
invariant algebraic curve.  In addition, we give an explicit
expression of a first integral for the system constructed from two
independent solutions of the second order ordinary differential
equation. This first integral is not, in general, a Liouvillian
function. The inverse integrating factor of the system (\ref{eq1})
which is associated to the aforementioned second order linear
differential equation, takes the form $V(x,y)=q(x)
(w'(x)-g(x,y)w(x))^2$ where $q(x)=A_2(x) \exp\left\{ \int
A_1(x)/A_2(x) \, dx\right\}$, $g(x,y)$ is a fixed rational
function and $w(x)$ is a non--null solution of the second order
linear differential equation. \par Moreover, in the work
\cite{GGG} we also consider first order linear differential
equations: $A_1(x) w'(x) + A_0(x) w(x) =0$ with polynomial
coefficients and analogous results are obtained. The inverse
integrating factor of the system (\ref{eq1}) which is associated
to this first order linear differential equation, takes the form
$V(x,y)=A_1(x) g(x,y)(w(x)-a(x,y))$ where $a(x,y)$ is a function
defined in terms of $A_1(x)$ and $A_0(x)$, $g(x,y)$ is fixed
rational functions and $w(x)$ is a non--null solution of the first
order linear differential equation. \par Hence, in the work
\cite{GGG}, we give families of systems which, by construction,
linearize, because its corresponding inverse integrating factor is
obtained in terms of the solutions of a linear differential
equation. Moreover, the given families are very general since they
come from any rational change of variables.  The present work
arises as a reciprocal of the work \cite{GGG}, since we look for
systems which can be linearized, in the sense of Definition
\ref{def1}.
\par The goal of this work is to demonstrate an algorithm to
detect when a system is integrable (either inside the Liouvillian
class or not) by means of a linearization process. That is, we
target to find systems whose integrability passes through the
solutions of a linear differential equation. We proceed by giving
and explaining several examples which illustrate this process.
\par The examples that we study suggest that the integrability by
linearization of a polynomial system (\ref{eq1}) reduces to solve
linear differential equations of order $2$ or it falls into the
Darboux theory of integrability. \par The studied examples also
make arise the following questions: when a system is linearizable
(in the sense of Definition \ref{def1}) with a linear differential
equation of order $1$, does it always exist a rational change of
variables which transforms the system to an orbital equation which
is linear? In the same way, we can also ask whether when a system
is linearizable (in the sense of Definition \ref{def1}) with a
linear differential equation of order $2$, does it always exist a
rational change of variables which transforms the system to an
orbital equation which is of Riccati type? \par For the families
of systems studied in \cite{GGG}, the answer to the previous two
questions is affirmative.
\newline

The question of linearizability has attracted many authors since
the transformation of an ordinary differential equation or a
partial differential equation of any order by means of several
differential--algebraic manipulations to a linear differential
equation, gives in general the solution of the first, nonlinear
problem. Moreover, ordinary differential equations which linearize
come naturally in some physical applications, see \cite{FS} and
the references therein. In the work \cite{FS}, the question of
which ordinary differential equations (of any order) linearize
upon differentiation is addressed and some sufficient conditions
on the form of the equation are given. However, these sufficient
conditions are very restrictive over the equations and only very
special particular equations can satisfy them. We only consider
ordinary differential equations of first order, that is systems of
the form (\ref{eq1}) and we study several differential--algebraic
manipulations so as to get a linear equation which characterizes
its integrability.
\newline

We use two different methods to exhibit that a system is
linearizable: {\em equivalence} and {\em compatibility}. Both
methods start in the same way. We consider a system (\ref{eq1})
and we think of it as polynomial in one variable, for instance
$y$. Then, we take a polynomial in the variable $y$ of a certain
fixed degree and with arbitrary coefficients, which are functions
of the variable $x$, $\sum_{i=0}^{M} h_{i}(x) y^{i}$, and we
impose it to be an inverse integrating factor of the corresponding
system (\ref{eq1}). This condition gives rise to a system of
linear differential equations on the coefficients $h_i(x)$. In
general, this system of linear differential equations is
overdetermined. Several conditions on system (\ref{eq1}) can make
this system compatible and the way to choose these conditions is
what distinguishes between both methods. \par The {\em
compatibility method} is the most gross: we consider the system of
linear differential equations with variables $h_{i}(x)$ and we
uncouple the variables by means of differentiation and resultants.
We end up with an algebraic--differential condition on system
(\ref{eq1}). Although this method gives all the possible choices
of system (\ref{eq1}) to have an inverse integrating factor of the
prescribed form, it is usually too overwhelming to be carried out.
\par On the other hand, we can consider the {\em equivalence method}. This method is wiser and consists on uncoupling the
system of linear differential equations with variables $h_{i}(x)$,
trying to avoid differentiation and resultants and only using
substitution at each step. We end up with a number of linear
differential equations of certain order and of only one variable,
say $h_0(x)$, and we make these equations equivalent, that is, we
impose them to be the same equation and/or to be identically null
for some of them. This way gives certain particular conditions on
the system (\ref{eq1}) which are, usually, easy to satisfy. \par
We note that the conditions given by the equivalence method are
also contained in the conditions given by the compatibility method
but their determination is much easier when equivalence is
involved. \par

The method of equivalence gives rise to one linear ordinary
differential equation of a certain order $\ell$ for one variable,
which can be any $h_i(x)$ in the expression $\sum_{i=0}^{M}
h_{i}(x) y^{i}$. Each solution of this ordinary differential
equation gives, by back substitution, an inverse integrating
factor for system (\ref{eq1}). We would like to know the possible
values of $\ell$, that is, we ask whether we can linearize systems
(\ref{eq1}) by means of a linear ordinary differential equation of
any order $\ell$, with $\ell \geq 0$. We aim to know the values of
the order $\ell$ corresponding to a linear differential equation
such that each of its solutions cannot be expressed as a
polynomial on the solutions of an equation of lower order. The
following result shows that the order of such a linear
differential equation is at most $2$.
\begin{theorem} We assume that system {\rm (\ref{eq1})} has an inverse
integrating factor of the form: \begin{equation} V(x,y) \, = \,
\sum_{i=0}^{M} h_i(x) \, y^i, \label{eqq} \end{equation} where $M$
is a nonnegative integer number and $h_i(x)$ are analytic
functions in $x$, $i=0,1,\ldots,M$. We assume that the functions
$h_i(x)$, for $i=1,2,\ldots,M$, are polynomials in $h_0(x)$ and
its derivatives, and that $h_0(x)$ satisfies a linear differential
equation of order $\ell$, with $\ell$ a nonnegative integer, whose
solutions cannot be algebraically expressed in terms of the
solutions of an equation of lower order. Then, $\ell \leq 2$.
\label{res1}
\end{theorem}
{\em Proof.} We know that a linear ordinary differential equation
of order $\ell$ has a fundamental set of solutions with cardinal
$\ell$. That is, there are $\ell$ linearly independent solutions
to the equation. Assume that $\ell \geq 3$ and let $V_1(x,y)$,
$V_2(x,y)$ and $V_3(x,y)$ be three inverse integrating factors
each one constructed by using one of these linearly independent
solutions through the expression (\ref{eqq}). The quotients of two
of them give first integrals of system (\ref{eq1}):
$H_1(x,y)=V_1(x,y)/V_3(x,y)$ and $H_2(x,y)=V_2(x,y)/V_3(x,y)$.
These two first integrals need to be functionally dependent since
any first integral of a planar differential system like
(\ref{eq1}) is a function of another one. We are going to show
that, in fact, $H_1$ and $H_2$ are algebraically dependent, that
is, there exists a polynomial with real coefficients $P(z_1,z_2)$
such that $P(H_1,H_2)\equiv 0$. We consider the level curves of
each $H_i$: $V_i(x,y)-c_iV_3(x,y)=0$, with $i=1,2$, which are two
polynomials in $y$ because each $V_i(x,y)$, $i=1,2,3$, is a
polynomial in $y$. Let us take the resultant of the polynomials
$V_1(x,y)-c_1V_3(x,y)$ and $V_2(x,y)-c_2V_3(x,y)$ with respect to
$y$ and we denote it by $R(c_1,c_2,x)$. We remark that this
resultant is a polynomial in $c_1$ and $c_2$ and we are going to
show that it factorizes as $R(c_1,c_2,x)=P(c_1,c_2)
R_0(c_1,c_2,x)$ where $P(c_1,c_2)$ and $R_0(c_1,c_2,x)$ are
polynomials in $c_1$ and $c_2$. This factorization of
$R(c_1,c_2,x)$ is deduced by the fact that each $y$--root of
$V_1(x,y)-c_1V_3(x,y)=0$, of a fixed $c_1$, needs to correspond to
a value of $c_2$ such that the whole $y$--root is contained in
$V_2(x,y)-c_2V_3(x,y)=0$. That is, fixed a $c_1$ and a $y$--root
of $V_1(x,y)-c_1V_3(x,y)=0$, there exists a value of $c_2$ such
that this $y$--root is completely contained in
$V_2(x,y)-c_2V_3(x,y)=0$. Therefore, we have encountered a
polynomial $P(c_1,c_2)$ which relates the two first integrals in
the desired way. \par Let us call $S_i(x)$ the solution of the
linear ordinary differential equation which gives the inverse
integrating factor $V_i(x,y)$, $i=1,2,3$. Since $P(H_1,H_2) \equiv
0$, we deduce that there exists a homogeneous polynomial with real
coefficients such that $p(S_1,S_2,S_3) \equiv 0$. To encounter
this polynomial $p$, we put $y=0$ in the expression of
$P(H_1,H_2)$ and we take a common denominator. The existence of
this polynomial $p$ implies that $S_3$ can be algebraically
expressed in terms of $S_1$ and $S_2$. We remark that any two
given functions $S_1$ and $S_2$ satisfy the linear homogeneous
ordinary differential equation of second order:
\[ \det \left|
\begin{array}{ccc} w''(x) & w'(x) & w(x) \\ S_1''(x) & S_1'(x) &
S_1(x) \\ S_2''(x) & S_2'(x) & S_2(x) \end{array} \right| \, = \,
0. \] Thus, the function $S_3$ is algebraically expressed in terms
of the solutions of an equation of second order. \par We conclude
that any inverse integrating factor which shows the
linearizability of a system (\ref{eq1}) through a linear
differential equation of order $\ell$ with $\ell \geq 3$, can be
expressed in a way that the linearizability of the system is given
through a linear differential equation of, at most, second order.
\bbox \newline

In order to make more precise the given notion of linearization of
a system (\ref{eq1}), we include several examples of this
phenomenon and we use the described methods of equivalence or
compatibility.

\section{Examples \label{sect2}}

\subsection{Automatically linearizable systems \label{sect21}}

In this section we describe examples of systems for which we find
an inverse integrating factor constructed by means of the
solutions of a linear differential equation, that is, we describe
examples of linearizability. We do not need to impose any
condition on the system to ensure its linearizability, that is,
the system of linear differential equations on the functions
$h_i(x)$ for which $V(x,y)=\sum_{i=0}^{M} h_{i}(x) y^{i}$ is an
inverse integrating factor is not overdetermined. \newline

{\bf Example 1.} Let us consider systems of the form (\ref{eq1})
with:
\begin{equation} P(x,y) \, = \, -y, \qquad Q(x,y) \, = \, \sum_{i=0}^{m} g_{2i}(x) \,
y^{2i}, \label{eq2} \end{equation} where $g_{2i}(x)$ are analytic
functions and $m$ is an integer number with $m \geq 0$. In this
section we only consider the case in which $m=2$.
\par We look for an inverse integrating factor $V(x,y)$ which is a polynomial in $y$ of the same degree as (\ref{eq2}) and of the
form: \[ V(x,y) \, = \, \sum_{i=0}^{m} h_{2i}(x) \, y^{2i}, \]
where $h_{2i}(x)$ are suitable functions which will satisfy linear
differential equations. Our goal is to impose such a function
$V(x,y)$ as an inverse integrating factor for system (\ref{eq2})
and to deduce the relations on the functions $g_{2i}(x)$ to
accomplish it. \newline

The case $m=0$ is easily integrable since the corresponding
orbital equation has separated variables. Let us explicit the
computations made when $m=2$. We have the system:
\begin{equation} \dot{x} \, = \, -y, \quad \dot{y}\, = \, g_0(x) + g_2(x) y^2 + g_4(x) y^4,
\label{eqk2} \end{equation} and we look for an inverse integrating
factor of the form $V(x,y) = h_0(x) \, + \, h_2(x) y^2 \, + \,
h_4(x) y^4$. We impose the following relation to be satisfied:
\[ P \, \frac{\partial V}{\partial x} \,  +\, Q \, \frac{\partial V}{\partial
y} \, = \, \div \cdot V, \] where $\div= 2 g_2(x) y + 4 g_4(x)
y^3$. The previous partial differential equation can be arranged
in powers of $y$ and we get that three relations among the
functions $h_{2i}(x)$ need to be satisfied, corresponding to the
powers $y^5$, $y^3$ and $y$. We remark that we have three
relations and three functions to satisfy them. This is an
exceptional case since the number of equations and the number of
variables coincide. These relations read for:
\begin{equation}
\begin{array}{lll}
\displaystyle h_4'(x) \, - \, 2 g_2(x)h_4(x) + 2 g_4(x)h_2(x) & =
& 0,  \vspace{0.2cm} \\ \displaystyle
h_2'(x) \, - \, 4 g_0(x)h_4(x) + 4 g_4(x)h_0(x) & = & 0,  \vspace{0.2cm} \\
\displaystyle h_0'(x) \, - \, 2 g_0(x)h_2(x) + 2 g_2(x)h_0(x) & =
& 0. \end{array} \label{eqx2}
\end{equation}
We can deduce the values of the functions $h_2(x)$ and $h_0(x)$,
for instance, from the first two equations (\ref{eqx2}). The third
equation in (\ref{eqx2}) gives a third order linear differential
equation for $h_4(x)$ which is: \begin{equation}
\begin{array}{l} \displaystyle g_4^2 h_4''' \, - \, 3 g_4 g_4'
h_4'' \, + \, (-4g_2^2
g_4^2+16g_0g_4^3-4g_4^2g_2'+4g_2g_4g_4'+3g_4'^2-g_4g_4'') h_4' \,
+ \, \vspace{0.3cm} \\ \displaystyle + 2\Big(
(4(g_4g_0'-g_0g_4')-2g_2g_2'-g_2'')g_4^2   \, + \, \vspace{0.2cm}
\\ \displaystyle \qquad \qquad \, + \,
(2g_2^2+3g_2')g_4g_4'-3g_2g_4'^2+g_2g_4g_4'' \Big) h_4  =  0,
\end{array} \label{eq3} \end{equation} where we avoid to write the
dependence on the variable $x$ to simplify the notation.
\par We have obtained a third order linear differential equation
whose solutions give rise to an inverse integrating factor for
system (\ref{eq2}) with $m=2$. We note that we did not need to
impose any restriction on the functions $g_{2i}(x)$, $i=0,1,2$, so
as to linearize the system. System (\ref{eq2}) with $m=2$, that is
system (\ref{eqk2}), is always linearizable and we have been able
to deduce this fact using our method. \par Although we have
obtained a third order linear differential equation for $h_4(x)$,
we are going to see that this function is related to a second
order linear ordinary differential equation, as the statement of
Theorem \ref{res1} establishes. Let us impose an inverse
integrating factor for system (\ref{eqk2}) of the form
$V(x,y)=(\tilde{h}_0(x) + \tilde{h}_2(x) y^2)^2$, where
$\tilde{h}_0(x)$ and $\tilde{h}_2(x)$ are suitable functions.
Repeating the same computations as before, we get that the
relation corresponding to the power $y^5$ in the equation of
inverse integrating factor gives that $\tilde{h}_0(x) \, = \,
(g_2(x)\tilde{h}_2(x)-\tilde{h}_2'(x))/(2g_4(x))$ and the other
two relations (corresponding to $y^3$ and $y^1$) are equal and
give the following second order ordinary differential equation for
$\tilde{h}_2(x)$: \begin{equation} g_4 \, \tilde{h}_2'' \, - \,
g_4'\, \tilde{h}_2' \, + \, (g_2g_4'-g_2'g_4+4g_0g_4^2-g_2^2g_4)
\, \tilde{h}_2 \, = \, 0. \label{eqss}
\end{equation}
The fact that the two relations, corresponding to $y^3$ and $y^1$,
are equal is not expected and this equality confirms that system
(\ref{eqk2}) is linearizable. \par Actually, straightforward
computations show that if we denote by $A(x)$ and $B(x)$ two
independent solutions of the second order equation (\ref{eqss}),
then a fundamental set of solutions of (\ref{eq3}) is: $A^2(x)$,
$A(x)B(x)$ and $B^2(x)$. Hence, the third order linear ordinary
differential equation (\ref{eq3}) is, in fact, reducible to a
second order ordinary differential equation. This is an example of
the result stated in Theorem \ref{res1}. We note that, in any
case, we have linearized system (\ref{eqk2}). The reduction of
order passes, in this case, through the change $h_4(x) \, = \,
(\tilde{h_2}(x))^2$ which transforms equation (\ref{eq3}) in a
nonlinear differential equation which is compatible with
(\ref{eqss}). \par We need to make precise that when we have
linearized system (\ref{eq2}) with $m=2$, that is system
(\ref{eqk2}), we have obtained a third order linear differential
equation (\ref{eq3}), but this equation was not necessary since an
inverse integrating factor can be obtained through a second order
linear ordinary differential equation, as Theorem \ref{res1}
states. \par We consider the algebraic change of variables $y
\mapsto z$ with $z=\tilde{h}_0(x) + \tilde{h}_2(x) y^2$, where
$\tilde{h}_i$, $i=0,2$ are the functions which define the
encountered inverse integrating factor $V(x,y)=(\tilde{h}_0(x) +
\tilde{h}_2(x) y^2)^2$. This change of variables applied to system
(\ref{eqk2}) gives the following orbital equation: \[
\frac{dz}{dx} \, = \, - \, \frac{z
\left(2zg_4(x)+\tilde{h}_2'(x)\right)}{\tilde{h}_2(x)} \, ,\]
which is of Riccati type. \par An algebraic change applied to
system (\ref{eqk2}) which gives an orbital equation of Riccati
type is not unique. We note that the algebraic change of variables
$y \mapsto u$ with $u=y^2$ applied to system (\ref{eqk2}) gives
the following orbital equation:
\begin{equation} \frac{du}{dx} \, = \, -2\, (g_0(x) \, + \, g_2(x)
\, u \, + \, g_4(x) \, u^2), \label{eqric}
\end{equation} which is also of Riccati type. The linearizable systems studied in the examples show
the existence of a rational change of variables which transforms
the system to an orbital equation of Riccati type, as a
counterpart. We do not look for changes of variables but for
linearizability.
\par It is well--known the equivalence between ordinary
differential equations of Riccati type and second order linear
differential equations. For instance, by the change $u \mapsto w$
with $dw/dx\, = \, 2 g_4(x) \, u(x) \, w$ we have that the Riccati
equation (\ref{eqric}) is equivalent to the following second order
linear differential equation for $w(x)$:
\begin{equation} g_4(x) \, w''(x) \, + \,
\left(2g_2(x)g_4(x)-g_4'(x)\right) w'(x) \, + \, 4 g_0(x) g_4^2(x)
\, w(x) \, = \, 0. \label{eqsec} \end{equation} \par This family
of systems also appears as a particular case of the systems
described in \cite{GGG} by means of the algebraic change of
variables $y \mapsto u$ with $u=y^2$, which leads it to the
orbital equation of Riccati type (\ref{eqric}). Following the
ideas described in \cite{GGG}, if we consider equation
(\ref{eqss}) and we perform the change $\tilde{h}_2(x)=\exp
\left\{ \int_{0}^{x} g_2(s) ds \right\} w(x)$, we obtain the
linear differential equation (\ref{eqsec}). We would like to
remark that the linearizability process does not pass through the
second order linear differential equation (\ref{eqsec}), which is
equivalent to the orbital equation (\ref{eqric}). The
linearizability of system (\ref{eqk2}) is concerned with the
second order linear differential equation (\ref{eqss}) whose
solutions define an inverse integrating factor for the system.
However, to impose that system (\ref{eqk2}) has an inverse
integrating factor constructed with the solutions of a linear
differential equation seems to imply that the orbital equation
associated to the system is equivalent to a linear differential
equation by means of an algebraic change of variables. All the
examples that we present in this work confirm this implication
although the involved linear differential equations are different
and come from different sources.

\subsection{Linearizability by equivalence \label{sect22}}

In this section we describe several examples of systems of the
form (\ref{eq1}) which, under certain restrictions, are
linearizable. The way we determine these restrictions is by using
the equivalence method. \newline

{\bf Example 2.} We consider system (\ref{eq2}) with $m \geq 3$
and we explain a process which encounters a linearizable
subfamily. Let us describe the computations for system (\ref{eq2})
with $m=3$. In fact, for $m=1$ and for any $m \geq 3$ the
discussion is analogous. The only cases which are different are
$m=0$ (separated variables) and $m=2$ (second order linear
differential equation) which have already been treated in Example
1. We consider the system:
\[ \dot{x} \, = \, -y, \quad
\dot{y}\, =\, g_0(x) + g_2(x) y^2 + g_4(x) y^4 + g_6(x) y^6, \]
and we take a function of the form $V(x,y)\, =\, h_0(x) \, + \,
h_2(x) y^2 \, + \, h_4(x) y^4 \, + \, h_6(x) y^6$. By imposing it
to be an inverse integrating factor and equating the same powers
of $y$, we get $5$ relations which need to be satisfied and which
correspond to the coefficients of $y$, $y^3$, $y^5$, $y^7$ and
$y^{9}$. From the coefficient of $y^{9}$ we compute $h_4(x)$ in
terms of $g_{2i}(x)$ and $h_6(x)$. In the same way, from the
coefficient of $y^7$ we compute $h_2(x)$ and from the coefficient
of $y^5$, we compute $h_0(x)$ in terms of $g_{2i}(x)$ and
$h_6(x)$. We are left with two linear differential equations of
second order for $h_6(x)$ which read for:
\[ \begin{array}{l}
\displaystyle   \frac{h_6''}{4g_6} \, + \,
\left(\frac{g_4^2}{3g_6^2} - \frac{g_2}{g_6} -
\frac{g_6'}{4g_6^2}\right) h_6' \, + \, \left[ \frac{1}{3g_6}
\left( \frac{g_4^2}{g_6} \right)' \, - \, \frac{g_4^2g_6'}{3g_6^3}
\, - \, \left(\frac{g_2}{g_6}\right)' \, \right] h_6 \, = \, 0,
\vspace{0.3cm} \\
\displaystyle \frac{g_4 \, h_6''}{12g_6} \, + \, \left(
\frac{1}{3} \left(\frac{g_4}{g_6}\right)' - \frac{g_4'}{12g_6} +
\frac{g_2g_4}{6g_6} - \frac{3g_0}{2} \right) h_6' \, +
\vspace{0.2cm} \\ \displaystyle \qquad  \quad + \, \left(
\frac{g_4g_6'^2}{2g_6^3}  - \frac{g_4'g_6'}{2g_6^2} +
\frac{g_2}{3} \left(\frac{g_4}{g_6}\right)' + \frac{g_0g_6'}{g_6}
- g_0' \right) h_6 \, = \, 0.
\end{array} \] We apply the equivalence method to these two linear
differential equations, that is, we impose the values of
$g_{2i}(x)$ to make them the same equation. Astonishingly, we only
need to impose the condition:
\[ g_0(x) \, = \, \frac{g_2(x)g_4(x)}{3g_6(x)} \, - \,
\frac{2}{27} \frac{g_4(x)^3}{g_6(x)^2} \, + \, \frac{1}{6} \left(
\frac{g_4(x)}{g_6(x)} \right)' , \] so as to get only one second
order linear differential equation for $h_6(x)$: \begin{equation}
3 g_6 h_6'' \, + \, (4g_4^2-12g_2g_6-3g_6')h_6' \, + \, \left(
8g_4g_4'-\frac{8g_4^2g_6'}{g_6} + 12 g_2 g_6' - 12 g_6 g_2'\right)
h_6 \, = \, 0. \label{eqss3} \end{equation} Hence, we have a
family of systems of the form (\ref{eq2}) which linearize. At this
point we have met our target since we have encountered a family of
systems whose integrability passes through the solution of a
second order linear differential equation.
\par We are going to describe another unexpected phenomenon which occurs in
this family. The second order linear equation for $h_6$ can be
reduced to a linear equation of order $1$. We remark that Theorem
\ref{res1} does not apply because equation (\ref{eqss3}) is of
order $2$. We consider system (\ref{eq2}) with $m=3$ and the
described value of the function $g_0(x)$. We impose a function of
the form $V(x,y)=(\tilde{h}_0(x)+\tilde{h}_2(x)y^2)^3$ to be an
inverse integrating factor. We get that
$\tilde{h}_0=g_4\tilde{h}_2/(3g_6)$ and we only obtain one linear
homogeneous differential equation of order $1$ for the function
$\tilde{h}_2(x)$ which is $9 g_6 \tilde{h}_2' \, + \,
4(g_4^2-3g_2g_6) \tilde{h}_2 \, =0$. Therefore, we have that
system (\ref{eq2}) with $m=3$ and the described value of $g_0(x)$
linearizes. We remark that the final equation for $\tilde{h}_2(x)$
in the case $m=3$ is of order $1$ whereas the final equation for
$\tilde{h}_2(x)$ in the case $m=2$ is of order $2$.
\par Moreover, this family of systems also appears as
a particular case of the work \cite{GGG}. We consider the
algebraic change of variables $y \mapsto z$ with $z=\tilde{h}_0(x)
+ \tilde{h}_2(x) y^2$, where $\tilde{h}_i$, $i=0,2$ are the
functions which define the encountered inverse integrating factor
$V(x,y)=(\tilde{h}_0(x) + \tilde{h}_2(x) y^2)^3$. This change of
variables applied the considered family of systems gives the
following orbital equation: \begin{equation} \frac{dz}{dx} \, = \,
 \frac{2}{9} \, z \left( -3 g_2(x) \, + \,
\frac{g_4(x)^2}{g_6(x)} \, - \, \frac{9 g_6(x)
z^2}{\tilde{h}_2(x)^2} \right), \label{eqber} \end{equation} which
is a differential equation of Bernoulli type. The change $z
\mapsto u$ with $z=1/\sqrt{u}$ transforms equation (\ref{eqber})
to a linear differential equation. We note that by this method we
impose a system to be linearized, that is, the other way round of
what we obtained in \cite{GGG}, where we started by the linear
differential equation and we deduced the corresponding system. As
far as the examples described in this work indicate, we have
observed that when a linearizable system (\ref{eq1}) (in the sense
of Definition \ref{def1}) admits an inverse integrating factor of
the form $V(x,y)=c(x,y)^n$, where $n$ in a positive integer and
$c(x,y)$ is a polynomial in $y$, the algebraic change $y \mapsto
z$ where $z=c(x,y)$ applied to the system gives an orbital
equation which is of Riccati or Bernoulli type.
\par The same computations can be done for the systems (\ref{eq2})
with $m=1$ or $m >3$, and giving certain values to $m-2$ of the
functions $g_{2i}(x)$, that is we have $3$ arbitrary functions
$g_{2i}(x)$ in the expression of $Q(x,y)$, an inverse integrating
factor of the form $V(x,y)=(\tilde{h}_0(x)+\tilde{h}_2(x)y^2)^m$
is exhibited, where $\tilde{h}_0(x)$ is expressed in terms of the
$g_{2i}(x)$ and $\tilde{h}_2(x)$, and the function
$\tilde{h}_2(x)$ is the solution of a linear homogeneous ordinary
differential equation of order $1$. In the particular case where
the linearizable system (\ref{eq2}) with $m \geq 3$ is a
polynomial differential system, we have that the encountered
inverse integrating factor is of Darboux type. We remark that in
the method of linearization we are not seeking for invariant
algebraic curves but we have obtained them by construction.
\par When we have
studied systems of the form (\ref{eq2}) with $m=2$ we have
exhibited a system in which the linearization is automatically
met. When we have described the systems of the form (\ref{eq2})
with $m=3$, we have given a family of linearizable systems in
which the linearization is met by the equivalence method. We also
remark that in case system (\ref{eq2}) is polynomial, we can
encounter non Liouvillian inverse integrating factors when $m=2$,
that is, when the system is of degree $4$ in $y$. When $m \geq 3$
in a polynomial system (\ref{eq2}), the linearization process
gives a Darboux inverse integrating factor and, thus, a
Liouvillian first integral. This fact implies that the encountered
second order linear differential equation (\ref{eqss3}) only has
elementary solutions. \newline

{\bf Example 3.} The following family of systems is quite more
general and we are going to use equivalence so as to encounter a
subfamily of linearizable systems. Let us consider the system:
\begin{equation} \begin{array}{lll}
\dot{x} & = & c_0(x) \, + \, c_2(x) y^2, \vspace{0.2cm} \\ \dot{y}
& = & d_0(x) + d_1(x) y + d_2(x) y^2 + d_3(x) y^3 + d_4(x) y^4 +
d_5(x) y^5 + d_6(x) y^6, \end{array} \label{eqg3} \end{equation}
where $c_i(x)$ and $d_i(x)$ are analytic functions with $c_2(x)
\cdot d_6(x) \not\equiv 0$. \par We start by imposing an inverse
integrating factor of the form \[ V(x,y) \, = \, \sum_{i=0}^{6}
h_i(x) \, y^i \, , \] which is a polynomial in $y$ of the same
degree as system (\ref{eqg3}). We substitute this expression of
$V$ in the partial differential equation that must be satisfied to
be an inverse integrating factor and we obtain that a certain
polynomial in $y$ of degree $10$ must be identically zero. From
the coefficients in $y$ of this equation of degrees from $10$ to
$5$ we deduce the values of $h_i(x)$ with $i=0,1,\ldots, 5$. The
rest of the coefficients give five linear differential equations
for $h_6(x)$. Two of them are of order $3$ and the rest are of
order $2$. We impose the three equations of order $2$ to be
identically null, which give the following conditions on system
(\ref{eqg3}):
\begin{equation}
\begin{array}{c} \displaystyle d_5(x) \, = \, 0, \quad d_4(x) \, = \,
\frac{6c_0(x)d_6(x)}{c_2(x)}, \quad d_2(x) \, = \,
\frac{9c_0(x)^2d_6(x)}{c_2(x)^2}, \vspace{0.2cm} \\ \displaystyle
d_1(x) \, = \, \frac{3c_0(x)d_3(x)}{c_2(x)} \, - \, c_2(x) \left(
\frac{c_0(x)}{c_2(x)} \right)' \, . \end{array} \label{eqcn}
\end{equation} Under these conditions we get that the three equations of
order $2$ are identically null and, surprisingly, the two
equations of third order are equal. The fact that under this small
number of restrictions, the involved linear differential equations
become equivalent is unexpected and confirms the hidden structure
of the linearizability process. This subfamily of systems
(\ref{eqg3}) is linearizable. \par We remark that this is not the
only way to proceed so as to get linearizability. We could also
have imposed the two equations of third order to be identically
null and then use equivalence with the second order equations. We
have only presented one of the possible cases that can be
encountered.
\par We know that the considered third order linear differential equation can be
reduced because its solutions cannot be functionally independent.
Let us consider an inverse integrating factor for system
(\ref{eqg3}) of the form $V(x,y) \, = \, (\tilde{h}_3(x) y^3 \, +
\, \tilde{h}_2(x) y^2 \, + \, \tilde{h}_1(x) y \, + \,
\tilde{h}_0(x))^2$. Analogous computations give that, under the
same conditions, the determination of $V(x,y)$ comes from the
solution of the following second order linear differential
equation in $\tilde{h}_3(x)$:
\[ \begin{array}{l}
\displaystyle \frac{2 c_2^2}{3 d_6} \, \tilde{h}_3'' \, - \,
\frac{2c_2^2d_6'}{3 d_6^2} \, \tilde{h}_3 ' \, + \, \left(
\frac{c_2'^2}{6 d_6} - \frac{c_2 c_2''}{3d_6} +
\frac{c_2c_2'd_6'}{3d_6^2} + c_2 \left(\frac{d_3}{d_6} \right)' -
\frac{3d_3^2}{2 d_6} + 6 d_0 \right) \tilde{h}_3 \, = \, 0.
\end{array} \]
We note that the change of variables $y \mapsto u$ with
\[ c_2(x)^2 u + 3 d_6(x)y(3c_0(x)+c_2(x)y^2)=0 \] transforms the
orbital equation associated to system (\ref{eqg3}), with the
values imposed in (\ref{eqcn}), to a Riccati equation. We have
that system (\ref{eqg3}), with the values described in
(\ref{eqcn}), is a particular case of the results given in
\cite{GGG}. As before, when we apply the linearization process we
do not look for systems which come from Riccati equations via a
change of variables, but we obtain such systems as a counterpart.

\subsection{Linearizability by compatibility \label{sect23}}

In this section we describe several examples of systems of the
form (\ref{eq1}) which, under certain restrictions, are
linearizable, in the sense of Definition \ref{def1}. We encounter
these restrictions by applying the compatibility method.
\newline

{\bf Example 4.} Let us consider the following system which
appears in the work \cite{CG}. The system
\begin{equation}
\dot{x} \, = \, y(-1+2\rho^2(x^2-y^2)), \quad \dot{y} \, = \,
x+\rho x^2+\rho y^2 + 4 \rho^2 x y^2, \label{eq40} \end{equation}
where $\rho \in \mathbb{R}$, has the inverse integrating factor
$V(x,y)=(x^2+y^2)^2(1 + 2 \rho x + \rho^2(x^2+y^2))$. We are going
to re--encounter this inverse integrating factor using
compatibility. This is the first example in which the
compatibility method appears. \par We consider system (\ref{eq40})
and we impose an inverse integrating factor of the form
$V(x,y)=h_0(x) +h_2(x)y^2 +h_4(x)y^4 +h_6(x)y^6$. We equate to
zero the coefficients of the powers of $y$ in the relation that
makes $V(x,y)$ an inverse integrating factor. We have five
equations corresponding to the coefficients of $y^j$ for
$j=1,3,5,7,9$ which read for:
\[ \begin{array}{l}
\displaystyle  2 x (1+ \rho x) h_2(x) \, = \, (1-2 \rho^2 x^2)
h_0'(x) + 2 \rho(1+6\rho x) h_0(x), \vspace{0.2cm} \\
\displaystyle 4 x (1+ \rho x) h_4(x) - 2 \rho^2 h_0'(x) \, = \,
(1-2 \rho^2 x^2) h_2'(x) +4 \rho^2 x h_2(x) , \vspace{0.2cm} \\
\displaystyle 6 x (1+ \rho x) h_6(x) - 2 \rho^2 h_2'(x)\, = \,
(1-2 \rho^2 x^2) h_4'(x) - 2 \rho (1+2\rho x) h_4(x), \vspace{0.2cm} \\
\displaystyle 2 \rho^2 h_4'(x)\, = \,
-(1-2 \rho^2 x^2) h_6'(x) + 4 \rho (1+3\rho x) h_6(x), \vspace{0.2cm} \\
\displaystyle 2 \rho^2 h_6'(x)\, = \, 0. \end{array} \] We remark
that we have five linear differential equation for the four
functions $h_0(x)$, $h_2(x)$, $h_4(x)$ and $h_6(x)$. This system
of linear differential equations is shown to be overdetermined. We
can proceed in several ways but we do not loss any generality in
following one of them. We can, for instance, take the equations
from the last to the first one in the order they have been
written. We solve them, leaving an arbitrary constant at each
step. We end up with several algebraic relations for these
constants which mark their value and from which we get the
previously described polynomial inverse integrating factor. The
compatibility of these relations is obtained by an adequate choice
of the constants of integration. Another way to study this system
of linear differential equations is to take them in the order they
have been written. From the first one we equate $h_2(x)$ and we
substitute in the rest of equations. From the second equation, we
equate $h_4(x)$ and we substitute in the rest of equations and
from the third equation we get $h_6(x)$. We end up with two linear
differential equations of fourth order for $h_0(x)$ and we make
them compatible. The compatibility process goes as follows: we
consider the two linear differential equations of fourth order for
$h_0(x)$ and we make a linear combination of them so as to get a
linear differential equation of third order. We derive this third
order linear differential equation and we combine it with one of
the previously considered linear differential equations of fourth
order for $h_0(x)$. We have two third order linear differential
equations for $h_0(x)$ and we combine them so as to get a linear
differential equation of second order. We derive it and we combine
with one of the previously considered equations of third order. We
have at this step two linear differential equations of second
order for $h_0(x)$, which we combine so as to get a first order
linear differential equation for $h_0(x)$. We derive it and we
obtain a second order linear differential equation for $h_0(x)$
which, combined with one of the previous gives rise to a first
order linear differential equation. The two first order linear
differential equations for $h_0(x)$ turn out to be the same. If
this was not  the case, we would combine them and we would obtain
a compatibility condition on the coefficients of the system. In
our case, we solve this first order linear differential equation
for $h_0(x)$. This value of $h_0(x)$ is the value which makes
compatible the two fourth order linear differential equations from
which we started the process. The only possible common solution of
these two fourth order linear differential equations is
$h_0(x)=x^4(1+\rho x)^2$ (modulus a multiplicative constant). We
observe that this $h_0(x)$ univocally determines the previously
described inverse integrating factor.
\par This example suggests that the integrability by linearization of a polynomial system
(\ref{eq1}) reduces to solve linear differential equations of
order $2$ or it falls into the Darboux theory of integrability. We
remark that any Darboux inverse integrating factor which is a
polynomial in the variable $y$ is encountered by our linearization
process: either by the equivalence method, either by the
compatibility method. \par We note that when applying the
linearization process to system (\ref{eq40}) we obtain a Darboux
inverse integrating factor and, thus, invariant algebraic curves
of the systems as a counterpart. \newline

{\bf Example 5.} In this example we are addressed the question
whether the family (\ref{eq40}) can be embedded in a linearizable
family and it can be seen as a particular case of linearization
with equivalence. \par We can think of system (\ref{eq40}) as a
particular case of the following family of systems:
\begin{equation}
\dot{x} \, = \, y \left( g_1(x) + g_2(x) y^2 \right), \quad
\dot{y} \, = \, g_3(x) + g_4(x) y^2, \label{eq4g}
\end{equation}
where the $g_i(x)$ are arbitrary functions. In case $g_2(x) \equiv
0$, we have that system (\ref{eq4g}) coincides with system
(\ref{eq2}) with $m=2$ after a time--rescaling, which has already
been studied in the first example. Therefore, we can assume,
without loss of generality, that $g_2(x)\equiv 1$. We note that
the family (\ref{eq4g}) is reversible by the change $(x,y,t)
\mapsto (x,-y,-t)$ and, thus, any ansatz of an inverse integrating
factor must be even in $y$. Let us consider an inverse integrating
factor of the form:
\[ V(x,y) \, = \, h_0(x) \, + \, h_2(x) y^2 \, + \, h_4(x) y^4 \,
+ \, h_6(x) y^6. \] The imposition for it to be an inverse
integrating factor gives a polynomial in $y$. We observe that this
relation is an odd polynomial in $y$ of degree $9$ and which is
odd. The functions $h_i(x)$ must vanish each one of the
coefficients of this polynomial in $y$. From the coefficient of
$y^9$ of this polynomial we deduce that $h_6(x) \, = \, k_6$, with
$k_6$ a constant value which we assume to be nonzero. From the
coefficient of $y$, we deduce the value of $h_2(x)$ and from the
coefficient of $y^3$, we deduce the value of $h_4(x)$. The
coefficients of $y^5$ and $y^7$ give rise to two linear
differential equations of order $3$ for $h_0(x)$.
\par When applying equivalence to these two equations, that is, imposing them to be the same equation, we deduce the following
conditions:
\[ g_3(x) \, = \, \frac{1}{4} \left(2g_1(x)g_4(x) - g_1(x)g_1'(x)\right),
\quad g_4(x) \, = \, - \, \frac{g_1'(x)}{2}. \] In this case, the
orbital equation associated to system (\ref{eq4g}) is of separated
variables. We observe that the linearization process for this
example leads to an ordinary differential equation with separated
variables.
\par Another way to study the possible linearizability of system
(\ref{eq4g}) is to impose that the two linear equations of third
order have a nonzero common solution, that is, to impose
compatibility. To perform the computations of compatibility for
two equations of third order carries long calculations and many
cases. If we apply compatibility to the two linear differential
equations of order $3$ for $h_0(x)$, we obtain several conditions.
The conditions obtained in the equivalence case are
re--encountered now. Moreover, we obtain two additional, and very
complicated, conditions on the functions $g_i(x)$ to have
compatibility. One of these two conditions is the one satisfied by
system (\ref{eq40}) as a particular case. We have seen that system
(\ref{eq40}) cannot be seen as a particular case of a linearizable
family of systems (\ref{eq4g}) by the equivalence method. \newline

{\bf Example 6.} We consider a planar differential system of the
form:
\begin{equation}
\dot{x} \, = \, g_0(x) + g_1(x) y + g_2(x) y^2, \quad \dot{y} \, =
\, g_3(x) y + g_4(x) y^2, \label{eq66} \end{equation} where
$g_i(x)$, $i=0,1,2,3,4$, are arbitrary functions. We remark that
this system has $y=0$ as invariant algebraic curve and we propose
an inverse integrating factor which contains this information. We
are going to give conditions on the functions $g_i(x)$,
$i=0,1,2,3,4$, such that the system has an inverse integrating
factor of the form:
\[ V(x,y) \, = \, h_1(x) y + h_2(x) y^2 + h_3(x) y^3, \] where
$h_i(x)$, $i=1,2,3$ are suitable functions. The imposition for $V$
to be an inverse integrating factor of system (\ref{eq66}) gives
rise to a polynomial in $y$ of degree $5$. From the coefficients
of $y^5$ and $y^1$ of this polynomial we deduce that:
\[ h_3(x) \, = \, k_2 \, g_2(x), \quad h_1(x) \, = \, k_0 \,
g_0(x), \] where $k_0$ and $k_2$ are arbitrary constants. We end
up with only three conditions which involve the functions
$g_i(x)$, $i=0,1,2,3,4$ and the constants $k_0$,$k_2$. Since the
function $h_2(x)$ is not concerned, we take $h_2(x) \equiv 0$ for
simplicity. The vanishing of this three conditions gives the
following planar differential system:
\begin{equation}
\begin{array}{lll}
\dot{x} & = & \displaystyle k_2 g_1(x)
(g_0(x)^2+g_0(x)g_1(x)y+g_1(x)^2y^2), \vspace{0.2cm} \\ \dot{y} &
= & \displaystyle y (k_0g_0(x) - k_2 g_0(x) -k_2 y
g_1(x))(g_0'(x)g_1(x)-g_0(x)g_1'(x)), \end{array} \label{eqf2}
\end{equation}
where $k_0,k_2$ are real constants and $g_0(x)$, $g_1(x)$ are
analytic functions. This system has the inverse integrating factor
$V(x,y)=y g_1(x)(k_0g_0(x)^2+k_2g_1(x)^2y^2)$. \par The following
rational change of variables $x \mapsto z$ with $z=1-k_2+y
g_1(x)/g_0(x)$ transforms the orbital equation associated to the
system to the following linear differential equation:
\[ \frac{dy}{dz} \, = \, \frac{y \,
z}{(1-k_2-z) \left(k_2 + k_1(1-k_2-z)^2\right)} . \]

\subsection{Inverse integrating factors of only one variable \label{sect24}}

{\bf Example 7.} In the following example we describe another
method of linearization, which consists on imposing conditions to
the system so as to obtain an inverse integrating factor depending
only on $y$. \par Let us consider the following system:
\begin{equation}
\dot{x} \, = \, y + c_1(x) + c_2(x) y^2, \quad \dot{y} \, = \,
c_3(x)+ c_4(x)y + c_5(x)y^2, \label{eqy1} \end{equation} where
$c_i(x)$, $i=1,2,3,4,5$ are arbitrary functions. Let us impose a
function $V=V_0(y)$ to be an inverse integrating factor and we
look for conditions on $c_i(x)$, $i=1,2,3,4,5$, to accomplish this
fact. The condition for $V_0(y)$ to be an inverse integrating
factor for (\ref{eqy1}) reads for:
\[ \begin{array}{l} \displaystyle \left( c_3(x)+ c_4(x)y + c_5(x)y^2 \right) V_0'(y) \, =
 \, \left( c_1'(x) +
c_4(x)+ 2 \, c_5(x)y +  c_2'(x) y^2 \right) V_0(y).
\end{array} \] Imposing that $c_j(x) \, = \, k_j \, c_1'(x)$ for $j=3,4,5$ and
$c_2(x) \, = \, k_1+ k_2 c_1(x)$ with $k_i$, $i=1,2,3,4,5$,
arbitrary constants, we get a linear differential equation for
$V_0(y)$. We change $c_1(x)$ to $c(x)$ for simplicity in notations
and we have that the system:
\begin{equation}
\dot{x} \, = \, c(x) + y + (k_1+k_2c(x))y^2, \quad \dot{y} \, = \,
c'(x)(k_3+k_4 y + k_5 y^2), \label{eqy2} \end{equation} linearizes
since it has an inverse integrating factor $V=V_0(y)$ which needs
to satisfy the following linear differential equation: $(k_3 +
k_4y+k_5y^2)V_0'(y) \, = \, (1+k_4 + 2 k_5y+k_2y^2)V_0(y)$. We
note that the function $V_0(y)$ is the exponential of the
primitive of a rational function, that is, it is of Darboux type.
\par

We are going to present a change of variables to system
(\ref{eqy2}) which transforms the corresponding orbital equation
to a linear differential equation. The following rational change
of variables $x \mapsto z$ with $z=c(x) + y + (k_1+k_2c(x))y^2$
transforms the orbital equation associated to system (\ref{eqy2})
into:
\[ \frac{dz}{dy} \, = \, \frac{1+2k_1y-k_2y^2}{1+k_2y^2} \, + \,
\frac{1+2k_2k_3y+2k_2y^2(1+k_4+k_5y)+k_2^2y^4}{(1+k_2y^2)(k_3+k_4y+k_5y^2)}
\, z, \] which is linear. Thus, we get that system (\ref{eqy2}) is
a particular case of the families described in \cite{GGG}.
\par
It is evident that linearizable systems of the form (\ref{eqy2})
of any degree in $y$, can be constructed in an analogous way. The
rational change of variables which would transform it to an
ordinary differential equation of linear type, and thus relate it
to the work \cite{GGG}, is $x \mapsto z$ with $z=P(x,y)$, where
$P(x,y)$ is the function defined by $\dot{x}=P(x,y)$. \par This
example suggests the following question: does it always exists an
algebraic change of variables that transforms a linearizable
system with a Darboux inverse integrating factor of only one
variable to a system whose orbital equation is linear?

\subsection{Inverse integrating factors of the form $V(x,y)=r(x)h(y)$. \label{sect25}}

{\bf Example 8.} The following examples of linearizable systems
are obtained imposing that the inverse integrating factor is a
product of two functions: one in the variable $x$ and the other in
the variable $y$.
\par Let us consider the system
\begin{equation}
\dot{x} \, = \, g_1(x) f_1(y) \, + \, g_2(x) f_2(y), \quad \dot{y}
\, = \, g_3(x) f_3(y) \, + \, g_4(x) f_4(y), \label{eqjj}
\end{equation}
where $g_i(x)$ and $f_i(y)$ are arbitrary functions, $i=1,2,3,4$.
We impose this system to have an inverse integrating factor of the
form $V(x,y)=r(x)h(y)$. When we substitute this expression in the
partial differential equation which defines an inverse integrating
factor, we get:
\[ \begin{array}{l} \displaystyle \left( f_3(y) + f_4(y) \, \frac{g_4(x)}{g_3(x)} \right)
h'(y) \, + \, \left[ f_4'(y) \, \frac{g_4(x)}{g_3(x)}  \, + \,
f_2(y) \left( \frac{g_2(x)r'(x)-g_2'(x)r(x)}{r(x)\, g_3(x)}
\right) \, + \right. \vspace{0.2cm} \\  \displaystyle  \left.
\qquad \qquad + \, f_1(y) \left( \frac{g_1(x)r'(x) - g_1'(x)
r(x)}{r(x)\, g_3(x)} \right) \right] h(y) \, = \, 0. \end{array}
\]
We take \[ g_3(x) \, = \, -\, \frac{g_1(x)}{k_0} \left(
\frac{g_2(x)}{g_1(x)} \right)', \  g_4(x) \, = \, -\,
\frac{g_1(x)}{k_0 k_1} \left( \frac{g_2(x)}{g_1(x)} \right)', \
r(x) \, = \, k_2 \, g_1(x), \] where $k_i$, $i=0,1,2$ are real
constants, and the previous relation reads for:
\[  \left( f_3(y) + \frac{f_4(y)}{k_1} \right)
h'(y) \, + \, \left( k_0 f_2(y)- f_3'(y) - \frac{f_4'(y)}{k_1}
\right) h(y) \, = \, 0, \] which is a linear ordinary differential
equation for $h(y)$. We obtain the following system:
\begin{equation}
\begin{array}{lll}
\dot{x} & = & \displaystyle k_0 \, g_1(x) \, (g_1(x) f_1(y) +
g_2(x) f_2(y)), \vspace{0.2cm} \\ \dot{y} & = &  \displaystyle
(k_1 f_3(y) + f_4(y))(g_1'(x)g_2(x)-g_1(x)g_2'(x))/k_1,
\end{array} \label{eqf1}
\end{equation}
where $k_0, k_1$ are real numbers and $g_1(x), g_2(x)$ and
$f_i(y)$, $i=1,2,3,4$, are arbitrary functions. This system has
the inverse integrating factor $V(x,y)=g_1(x)^2 h(x)$, where
$h(x)$ satisfies the aforementioned linear differential equation
of order $1$.
\par The following change of variables $x \mapsto z$ with
$z=g_2(x)/g_1(x)$ transforms the orbital equation associated to
system (\ref{eqf1}) into the following linear differential
equation:
\[ \frac{dz}{dy} \, = \, - \, k_0 \, \frac{f_1(y) \, + \, z \,
f_2(y)}{k_1 f_3(y) + f_4(y)}. \] As in the previous example we are
addressed the question of the existence of an algebraic change of
variables that transforms a linearizable system with a Darboux
inverse integrating factor of the form $V(x,y)=r(x)h(y)$ to a
system whose orbital equation is linear.

\section{More general inverse integrating factors}

In the first sections of this work we have provided several
examples to exhibit that certain families of systems (\ref{eq1}),
which are polynomial in the variable $y$, have an inverse
integrating factor which is a polynomial in $y$. In this section
we treat other expressions of an inverse integrating factor which
contain the polynomials in $y$ as a subclass. This generalization
is done in two steps. We first introduce a real parameter $\alpha$
which allows to study polynomial inverse integrating factors which
are a polynomial in $y$ up to a real power. This real parameter
does not involve any change in the linearization process described
so forth. The second step is to consider inverse integrating
factors which are a power series in $y$ and we end up with a
numerable set of linear differential--difference equations.
\par
We remark that the linearization process can also be carried out
by imposing a function $V_\alpha(x,y)$ of class $\mathcal{C}^1$ in
some open set $\mathcal{U}$ of $\mathbb{R}^2$, non locally null
and which satisfies the following partial differential equation:
\begin{equation} P \, \frac{\partial V_\alpha}{\partial x}\, +\,
Q\, \frac{\partial V_\alpha}{\partial y}\, = \, \alpha \left(
\frac{\partial P}{\partial x} + \frac{\partial Q}{\partial y}
\right) \, V_\alpha , \label{eqva} \end{equation} where $\alpha$
is a real number. In the case $\alpha=1$ we recover the method of
linearization described in the previous sections and this real
free parameter $\alpha$ gives a generalization of the
linearization process which can lead to wider families of
linearizable systems. The knowledge of a function $V_\alpha(x,y)$
satisfying this partial differential equations gives that $V\, =
\, V_\alpha^{1/\alpha}$ is an inverse integrating factor of system
(\ref{eq1}). In the following, we do not impose a system to have
an inverse integrating factor but to have a function
$V_\alpha(x,y)$ which satisfies the partial differential equation
(\ref{eqva}) where $\alpha$ is a real parameter.
\par In order to illustrate this linearization process  we
describe an example, where the expression of the function
$V_\alpha(x,y)$ is not a polynomial in $y$ but a power series in
$y$, that is,
\[ V_\alpha(x,y) \, = \, \sum_{n = 0}^{\infty} v_n(x) \, y^n. \]

Let us consider a quadratic polynomial differential system:
\begin{equation} \begin{array}{lll}
\dot{x} & = & \displaystyle a_{00} + a_{10} x + a_{01} y +
a_{20}x^2 + a_{11}xy + a_{02} y^2, \vspace{0.2cm} \\ \dot{y} & = &
\displaystyle b_{00} + b_{10} x + b_{01} y + b_{20}x^2 + b_{11}xy
+ b_{02} y^2,
\end{array} \label{eqqua}
\end{equation}
where $a_{ij}$ and $b_{ij}$ are real numbers. We impose that
$V_\alpha(x,y) \, = \, \sum_{n \geq 0} v_n(x) \, y^n$, satisfies
the corresponding partial differential equation (\ref{eqva}). We
fix a natural number $n$, we equate the coefficients of $y^n$ in
the development of (\ref{eqva}) and we get the following linear
differential--difference equation for $v_n(x)$. As in the rest of
the work, $v_n'(x)$ means the derivative of $v_n(x)$ with respect
to $x$.
\begin{equation}
\begin{array}{l}
\displaystyle \left( a_{00} + a_{10} x + a_{20} x^2 \right)
v_n'(x) \, + \, \left( a_{01} + a_{11} x \right) v_{n-1}'(x) \, +
\, a_{02} \, v_{n-2}'(x)  \vspace{0.2cm} \\ \displaystyle \, + \,
(n+1) \left( b_{00} + b_{10}x + b_{20} x^2 \right) v_{n+1}(x) \, +
\, \left[ (n-\alpha)(b_{01} + b_{11} x) \right. \vspace{0.2cm} \\
\left. \displaystyle \, - \, \alpha(a_{10} + 2 a_{20} x)
\right]v_{n}(x) \, + \, \left( (n-2\alpha -1) b_{02}\, - \, \alpha
a_{11} \right) v_{n-1}(x) \, = \, 0.
\end{array} \label{eqdd}
\end{equation} \newline

We include several examples of planar systems whose integrability
can be determined with a function $V_{\alpha}(x,y)$ which is a
power series in $y$. \newline

{\bf Example 9.} We are going to take $v_n(x)$ in recurrence
(\ref{eqdd}), of the form $v_n(x)\, = \, q(x) \, p_n(x)$ where
$q(x)$ is a suitable function in $x$ (usually of Darboux type) and
$p_n(x)$ is a polynomial in $x$. In particular, in this example we
impose a function of the form $v_{n}(x)\, = \, q(x) \, \varphi_n
\, H_n(x)$ to be a solution of equation (\ref{eqdd}), where $q(x)$
is a suitable function, $\varphi_n$ is a suitable sequence of real
numbers and $H_n(x)$ is the Hermite orthogonal polynomial of
degree $n$. For further information about orthogonal polynomials
and the identities they satisfy, see for instance
\cite{Abramowitz}. The Hermite orthogonal polynomials satisfy the
following two identities: \begin{equation} \begin{array}{l}
\displaystyle H_{n+1}(x) \, = \, 2 \, x \, H_n(x) \, - \, 2\, n\,
H_{n-1}(x), \qquad \displaystyle  H_n'(x) \, = \, 2 \, n \,
H_{n-1}(x).
\end{array} \label{idher} \end{equation} We substitute the expression $v_{n}(x)\, = \, q(x) \,
\varphi_n \, H_n(x)$ in (\ref{eqdd}) and we use the previous
identities to simplify it. We have that the following relation
must be fulfilled:
\[ a_{02} q(x) \varphi_{n-2} H_{n-2}'(x) \, + \, A_1(n,x) \,
H_{n-1}(x) \, + \, A_2(n,x) H_{n-2}(x) \, = \, 0, \] where the
$A_i(n,x)$ are expressions involving the sequence $\varphi_n$, the
function $q(x)$ and the parameters $a_{ij}$, $b_{ij}$ and
$\alpha$. Since we have already used all the identities of the
Hermite polynomials, we need that the coefficients of
$H_{n-2}'(x)$, $H_{n-1}(x)$ and $H_{n-2}(x)$ independently vanish.
We impose that $a_{02}=0$ and the following values make that
$A_i(n,x)=0$ for $i=1,2$:
\[ \begin{array}{c} \displaystyle \varphi_{n} \, = \, \frac{1}{n!}
\left(-\frac{a_{11}}{a_{20}}\right)^{n-1}, \ \ a_{00}=0, \ \
a_{10}=\frac{a_{01}a_{20}}{a_{11}}, \vspace{0.2cm} \\
\displaystyle  b_{01}=-\frac{a_{01}a_{20}}{a_{11}}, \ \ b_{11}= -
a_{20}, \ \ b_{02}=0, \ \ \alpha =-1, \end{array} \] \[
\displaystyle q(x) \, = \, \exp\left\{
\frac{2(a_{11}b_{10}-a_{01}b_{20})x+a_{11}b_{20}x^2}{a_{20}^2}\right\}
\left(a_{01}+a_{11}x\right)^{\mu},
\]
with $\ \displaystyle \mu \, = \, 2\,
\frac{a_{01}^2b_{20}}{a_{11}a_{20}^2} + 2 \,
\frac{(a_{11}b_{00}-a_{01}b_{10})}{a_{20}^2} -1$. \newline

We rename the free parameters by $a_{01}=a_0a_{11}$,
$a_{20}=-a_2a_{11}$, $b_{00}=b_0a_{11}$, $b_{10}=b_1a_{11}$,
$b_{20}=b_{2}a_{11}$ and $a_{11}$ is taken to be $a_{11}=1$. We
obtain that the quadratic system: \begin{equation} \dot{x} \, = \,
(a_0 + x)(y-a_2 x), \ \dot{y} \, = \, b_0 + b_1 x + b_2 x^2 + a_2
(a_0+x)y, \label{eqher} \end{equation} has the following
expression $V_{\alpha}(x,y)$ which satisfies the partial
differential equation (\ref{eqva}) with $\alpha=-1$:
\[ V_\alpha(x,y) \, = \, q(x) \, \sum_{n \geq 0} \varphi_n \,
H_n(x)\, y^n \, = \,  q(x)\, a_2 \, \sum_{n \geq 0} \frac{1}{n!}
\, H_n(x) \, \left(\frac{y}{a_2}\right)^n. \] The Hermite
polynomials have the following generating function:
\[ \exp\left\{ 2 x y - y^2 \right\} \, = \, \sum_{n \geq 0}
\frac{1}{n!} \,  H_{n}(x) \, y^n, \] which also appears in the
book \cite{Abramowitz}. This identity allows us to identify the
power series given by $V_{\alpha}(x,y)$ and we obtain the
following inverse integrating factor $V=V_\alpha^{-1}$ for system
(\ref{eqher})
\[ V_\alpha(x,y) \, = \, \exp\left\{ \frac{y^2-2a_2xy-b_2x^2+2(a_0b_2-b_1)x}{a_2^2}\right\} \,
 (a_0+x)^{\frac{-a_2^2-2b_0+2a_0b_1-2a_0^2b_2}{a_2^2}} , \]
which is of Darboux type and which is not a polynomial in any of
the variables $x$ or $y$. \par In this example, we have solved
recurrence (\ref{eqdd}) by imposing it to be compatible with the
identities (\ref{idher}). Using the solution of the recurrence, we
have encountered the inverse integrating factor of system
(\ref{eqher}) and we have, therefore, integrated the system.
\newline

{\bf Example 10.} The system $\  \dot{x} \, = \, 1-x^2, \ \dot{y}
\, = \, y(x-y),  $ has an inverse integrating factor of the form
$V=V_\alpha^{1/\alpha}$ with $\alpha=-1/2$ and where:
\[ V_\alpha(x,y) \, = \, \frac{1}{(1-x^2)^{1/4}} \, \frac{1-xy}{1-2xy+y^2}. \] This function $V_\alpha(x,y)$ has
been determined using that it satisfies;
\[ V_\alpha(x,y) \, = \, q(x) \, \sum_{n=0}^{\infty}
T_n(x) \, y^n, \] where $q(x)=(1-x^2)^{-1/4}$ and $T_n(x)$ is the
Chebyshev polynomial of first kind and of degree $n$. This choice
makes the recurrence equation (\ref{eqdd}) to be satisfied in this
case. As stated in the book \cite{Abramowitz}, the Chebyshev
polynomials of first kind satisfy the following identities:
\[ T_{n+1}(x) \, = \, 2\, x \, T_n(x) - T_{n-1}(x), \quad T_n'(x) \, =
\, 2nT_{n-1}(x) + \frac{n}{n-2} \, T_{n-2}'(x), \] \[
\frac{1-xy}{1-2xy+y^2} \, = \, \sum_{n=0}^{\infty} T_n(x) \, y^n.
\] \newline

{\bf Example 11.}  The system $\ \dot{x} \, = \, 1-x^2, \ \dot{y}
\, = \, 1-xy, $ has an inverse integrating factor of the form:
\[ V(x,y) \, = \, \frac{(1-x^2)^2}{\sqrt{1-2xy+y^2}}, \] which has
been determined using that it satisfies;
\[ V(x,y) \, = \, q(x) \, \sum_{n=0}^{\infty}
P_n(x) \, y^n, \] where $q(x)=(1-x^2)^{2}$ and $P_n(x)$ is the
Legendre polynomial of degree $n$. As stated in the book
\cite{Abramowitz}, the Legendre polynomials satisfy the following
identities:
\[ \begin{array}{c} \displaystyle (n+1) \, P_{n+1}(x) \, = \, (2n+1) \, x \, P_n(x) \, -\,  n \, P_{n-1}(x), \vspace{0.2cm} \\  \displaystyle  P_n'(x) \, =
\, \frac{n}{x^2-1} \, \left( x P_{n}(x) - P_{n-1}(x)\right),
\vspace{0.2cm} \\ \displaystyle  \frac{1}{\sqrt{1-2xy+y^2}} \, =
\, \sum_{n=0}^{\infty} P_n(x) \, y^n. \end{array}
\] \newline

{\bf Example 12.}  The system \begin{equation} \begin{array}{lll}
\dot{x} & = & \displaystyle
1+(2a^2-7)y^2+6y^4-2axy(1+y^2)+2x^2y^2,
\vspace{0.3cm} \\
\dot{y} & = & \displaystyle 2y^2(1-y^2)(a-xy),
\end{array} \label{eqhera} \end{equation} where $a$ is a real parameter has an inverse integrating factor of the form:
\[ V(x,y) \, = \, \frac{1}{\sqrt{1-y^2}} \, \exp\left\{ \frac{2axy-(a^2+x^2)y^2}{1-y^2} \right\}. \]
This function $V(x,y)$ has been determined using that it
satisfies;
\[ V(x,y) \, = \, \sum_{n=0}^{\infty}
\frac{H_n(a) \, H_n(x)}{n!\, 2^n} \, y^n, \] where $H_n(x)$ is the
Hermite polynomial of degree $n$. The identity \[
\frac{1}{\sqrt{1-y^2}} \, \exp\left\{
\frac{2axy-(a^2+x^2)y^2}{1-y^2} \right\} \, = \,
\sum_{n=0}^{\infty} \frac{H_n(a) \, H_n(x)}{n!\, 2^n} \, y^n, \]
is called {\em Mehler's--Hermite polynomial formula}, see for
instance \cite{special}.
\newline

These examples suggest that there is a connection between some
Darboux inverse integrating factors and orthogonal polynomials.

\vspace{0.5cm}

{\bf Addresses and e-mails:} \\
$^{\ (1)}$ Lab. de Math\'ematiques
et Physique Th\'eorique. CNRS UMR 6083. \\ Facult\'e des Sciences
et Techniques. Universit\'e de Tours. \\ Parc de Grandmont, 37200
Tours, FRANCE.
\\ {\rm E-mail:} {\tt Hector.Giacomini@lmpt.univ-tours.fr}
\vspace{0.2cm} \\
$^{\ (2)}$ Departament de Matem\`atica. Universitat de Lleida. \\
Avda. Jaume II, 69. 25001 Lleida, SPAIN. \\ {\rm E--mails:} {\tt
gine@matematica.udl.cat}, {\tt mtgrau@matematica.udl.cat}

\end{document}